\documentclass[paper=a4,english,fontsize=10pt,parskip=half,abstract=true]{scrartcl}
\usepackage{babel}
\usepackage[utf8]{inputenc}
\usepackage[T1]{fontenc}
\usepackage[left=20mm,right=20mm,top=30mm,bottom=30mm]{geometry}
\usepackage{amsmath}
\usepackage{amsthm}
\usepackage{amssymb}
\usepackage{mathtools}
\mathtoolsset{centercolon} 
\usepackage[bookmarks=true,
            pdftitle={An answer to a question of Bonnafe},
            pdfauthor={Benjamin Sambale},
            pdfkeywords={},
            pdfstartview={FitH}]{hyperref}

\numberwithin{equation}{section}

\allowdisplaybreaks[1]

\renewcommand{\phi}{\varphi}

\newcommand{\N}{\mathrm{N}}

\newcommand{\ZZ}{\mathbb{Z}}

\newcommand{\FF}{\mathbb{F}}

\newcommand{\Irr}{\operatorname{Irr}}

\title{An answer to a question of Bonnafé}
\author{Benjamin Sambale\footnote{Institut für Algebra, Zahlentheorie und Diskrete Mathematik, Leibniz Universität Hannover, Welfengarten 1, 30167 Hannover, Germany,
\href{mailto:sambale@math.uni-hannover.de}{sambale@math.uni-hannover.de}}}
\date{\today}

\begin{document}
\frenchspacing
\maketitle
\begin{abstract}\noindent
We give a negative answer to a question of Bonnafé on the Loewy length of a character ring of a finite group.
\end{abstract}
%

A question frequently addressed in group theory asks which properties of a finite group $G$ can be determined by its ordinary complex characters. The Grothendieck group of these characters becomes a ring $\ZZ\Irr(G)$, called the \emph{character ring}, with respect to tensor products. To make things more interesting consider a prime divisor $p$ of $|G|$ and a complete discrete valuation ring $\mathcal{O}$ such that $F=\mathcal{O}/J(\mathcal{O})$ is an algebraically closed field of characteristic $p$. We may take scalar extensions $\mathcal{O}\Irr(G)=\mathcal{O}\otimes_{\ZZ}\ZZ\Irr(G)$ and reduced modulo $J(\mathcal{O})$ to obtain $F\Irr(G)$. Bonnafé~\cite{BonnafeCR} has shown that $F\Irr(G)$ decomposes into local algebras, called blocks, and those are parametrized by the $p$-regular conjugacy classes of $G$. Recall that the $p'$-\emph{section} of a $p'$-element $g\in G$ consists of all $h\in G$ such that the $p'$-factor of $h$ is conjugate to $g$ in $G$. Let $e_S$ be the characteristic function of the $p'$-section $S$ of $g$ (i.\,e. $e_S(h)=1$ if $h\in S$ and $0$ otherwise). It can be shown that $e_S\in\mathcal{O}\Irr(G)$ and we may consider $e_S$ as an idempotent of $F\Irr(G)$. In fact, $e_S$ is a primitive idempotent, so every block of $F\Irr(G)$ has the form $B_S=F\Irr(G)e_S$ for some $p'$-section $S$. The block $B_p=B_p(G)$ corresponding to the $p'$-section $S=G_p$ of all $p$-elements is called the \emph{principal} block by Bonnafé. It should be mentioned that the more familiar block decomposition of $FG$ has little to do with the decomposition of $F\Irr(G)$ (not even the number of blocks is the same). 

Bonnafé discovered a remarkable similarity between $B_p(G)$ and $B_p(\N_G(P))$ where $P$ is a Sylow $p$-subgroup of $G$. 
He showed that both algebras have the same Loewy length whenever $P$ is abelian. This can perhaps be seen as a variation of McKay's, Alperin's or Broué's conjecture. After checking many more cases, Bonnafé asked at the end of his paper whether $B_p(G)$ and $B_p(\N_G(P))$ always have the same Loewy length.

The aim of this short note is to report on the following counterexample to Bonnafé's question: Let $p=2$ and let $G$ be the solvable group 
\[\mathtt{SmallGroup}(768,1085354)\cong C_8^2\rtimes C_3\rtimes C_4\] 
of order $768=2^8\cdot 3$. Here $B_p(G)$ has Loewy length $5$ while $B_p(\N_G(P))$ has Loewy length $6$.
The computation was performed with GAP~\cite{GAP48} using the following strategy: The character ring $\ZZ\Irr(G)$ can be realized as a structure constant algebra with respect to the canonical basis $\Irr(G)$. The structure constants are the scalar products $[\chi\psi,\phi]$ where $\chi,\psi,\phi\in\Irr(G)$. We only need to reduce these integers modulo $p$ to obtain $\FF_p\Irr(G)$. Similarly, the coefficient of $e_{G_p}$ with respect to $\chi\in\Irr(G)$ is $\frac{1}{|G|}\sum_{g\in G_p}\chi(g)$; a rational number by elementary Galois theory. Since $e_{G_p}\in\mathcal{O}\Irr(G)$, these coefficients can again be reduced modulo $p$. In this way we construct the subalgebra $\FF_p\Irr(G)e_{G_p}$. Now the Loewy length does not change under the scalar extension $B_p(G)=F\otimes_{\FF_p}\FF_p\Irr(G)e_{G_p}$ since $J(B_p(G))=F\otimes_{\FF_p}J(\FF_p\Irr(G)e_{G_p})$. 
Finally, the computation of the Loewy length was done with the command \texttt{ProductSpace} in GAP. According to the \texttt{SmallGroupsInformation} command, the first 1085323 groups of order $768$ in the small groups library have a normal Sylow $2$-subgroup or a normal Sylow $3$-subgroup. In those cases Bonnafé's question has an affirmative answer by \cite[Proposition~4.7]{BonnafeCR} (in part (d) of that proposition $H$ should be replaced by $G/N$). 

\section*{Acknowledgment}
I thank Gabriel Navarro for drawing my attention to Bonnafé's paper. He and Radha Kessar insisted that I put my findings on arXiv.
I am supported by the German Research Foundation (\mbox{SA 2864/1-2} and \mbox{SA 2864/3-1}).


\begin{thebibliography}{1}

\bibitem{BonnafeCR}
C. Bonnaf\'{e}, \textit{On the character ring of a finite group}, in: Alg\`ebre
  et th\'{e}orie des nombres. {A}nn\'{e}es 2003--2006, 5--23, Publ. Math. Univ.
  Franche-Comt\'{e} Besan\c{c}on Alg\`ebr. Theor. Nr., Lab. Math. Besan\c{c}on,
  Besan\c{c}on, 2006.

\bibitem{GAP48}
The GAP~Group, \textit{GAP -- Groups, Algorithms, and Programming, Version
  4.11.0}; 2020, (\url{http://www.gap-system.org}).

\end{thebibliography}
\end{document}